\documentclass[a4paper, 12pt, leqno]{amsart}
\usepackage{amsmath}
\usepackage{amscd}
\usepackage{amssymb}
\usepackage{amsthm}
\theoremstyle{definition}

\theoremstyle{remark}

\numberwithin{equation}{subsection}

\begin{document}

\hsize=14truecm
\title[Maximal and Primitive Elements]
 {  Maximal    and  Primitive        Elements\\
  in Baby Verma Modules for  Type $B_2$
}

\author{ Nanhua Xi }

\address{Institute of Mathematics, Chinese Academy of Sciences, Beijing 100080,
China}

\email{nanhua@math.ac.cn}


\thanks{The author was  supported in part by the National Natural Science Foundation
of China.}


\dedicatory{}

\maketitle



\baselineskip=18pt

\hoffset=-2truecm
\medskip\bigskip

\medskip
\bigskip
The purpose of this paper is to find maximal and primitive elements
of baby Verma modules  for a quantum group of type $B_2$.  As a
consequence the composition factors of the baby Verma modules are
determined. Similar approach can be used to find find maximal and
primitive elements of Weyl modules for type $B_2$. In principle the
results can be used to determine the module structure of a baby
Verma module, but the calculations are rather involved, much more
complicated than the case of type $A_2$.

For type $A_2$, submodule structure of a Weyl module has been
determined
  in [DS1, I, K] and by Cline
(unpublished).  For type $B_2$, the  socle series of Weyl modules
was determined in [DS2].  In [X2] we determine the maximal and
primitive elements in Weyl modules for type $A_2$, so that the  Weyl
modules are understood more explicitly. This paper is a sequent work
of [X2], but less complete, since submodule  structure of a baby
Verma module is not determined.  In this paper we only work with
quantized enveloping algebras at roots of 1 (Lusztig version). For
hyperalgebras the approach is completely similar, actually simpler.

The contents of the paper are as follows. In section 1 we recall
some definitions and results about maximal and primitive elements.
In section 2 we recall some facts about a quantized enveloping
algebra of type $B_2$. In section 3 we determine the maximal and
primitive elements in a Verma module of the (slightly enlarged)
Frobenius kernel of type $B_2$. In section 4 we indicate that the
maximal and primitive elements in a Weyl module for type $B_2$ can
be worked out similarly, but we omit the results. To avoid
complicated expressions and for simplicity we assume that the order
of the involved root of 1 is odd and greater than 3 and we only work
with some special weights.  The approach for general cases is
completely similar.

\bigskip

\section{    Maximal    and  Primitive Elements}

In this section we fix notation and recall the definition and some
results  for maximal and primitive          elements. We refer to
[L1-4, X1-2] for additional information.

\def\Uxi{U_\xi}
\def\qxi{\bold Q(\xi)}
\def\bu{\bold u}
\def\bzn{\bold Z^n}
\def\l{\lambda}
\def\ll{_\l}
\def\tbu{\tilde{\bu}}
\def\k{k_i}
\def\kk{k_i^{-1}}
\def\a{\alpha}
\def\b{\beta}
\def\uxin{\bold u^-}
\noindent{\bf 1.1.} Let $\Uxi$ be a quantized enveloping algebra
 (over $\qxi$)
at a root $\xi$ of 1 (Lusztig version). We assume that the rank of
the associated Cartan matrix is $n$ and the order of $\xi\geq 3$. As
usual, the generators of $\Uxi$ are denoted by $e_i^{(a)},
f_i^{(a)},\k,\kk$, etc. Let $\bu$ be the Frobenius kernel and $\tbu$
the subalgebra of $\Uxi$ generated by all elements in $\bu$ and in
the zero part of $\Uxi$. For $\l\in\bzn$ and a $\Uxi$-module (or
$\tbu$-module $M$) we denote by $M\ll$ the $\l$-weight space of $M$.
A nonzero element in $M\ll$ will be called a vector of weight $\l$
or a weight vector. Let $m$ be a weight vector of a $\Uxi$-module
(resp. $\tbu$-module) $M$. We call $m$ {\bf maximal} if
$e_i^{(a)}m=0$ for all $i$ and $a\geq 1$ (resp. $e_\a m=0$ for all
root vectors $e_\a$ in the positive part of $\tbu$). We call $m$ a
{\bf primitive
     element}  if there exist two
submodules $M_2\subset M_1$ of $M$ such that $m\in M_1$ and
the image in $M_1/M_2$ of $m$ is maximal.     Obviously, maximal
    elements are  primitive.     We
have (see [X2]:

\medskip

\noindent (a) Let $m\in M$ be a weight vector and let $P_1$ be the
 submodule of $M$ generated by $ m$. Then $m$ is  primitive          if and only if the image in $P_1/P_2$ of $m$ is
 maximal     for some proper submodule $P_2$ of $P_1$.

We shall write $L(\l)$ (resp. $\tilde L(\l)$)
 for an irreducible $\Uxi$-module (resp. $\tbu$-module) of highest weight $\l$.
\def\tll{\tilde L(\l)}

\noindent (b) If $m$ is a  primitive          element of weight
$\l$, then $L(\l)$ (or $\tilde L(\l)$) is a composition factor of
$M$ (depending on $M$ is a $\Uxi$-module or a $\tbu$-module).

\noindent (c) Let $M$ and $N$ be modules and $\phi: M\to N$ a homomorphism.
Let $m$ be a weight vector in $M$. If $\phi(m)$ is a  primitive
element of $N$, then $m$ is a  primitive          element of $M$.

\noindent (d) Let $M,N, \phi, m$ be as in (c) and assume $\phi(m)\ne 0$ . If $m$ is  a primitive
element of $M$,
 then either $\phi(m)$ is a  primitive          element of $N$ or $\phi(P_1)=\phi(P_2)$, where
$P_1$ is the submodule of $M$ generated by $m$ and $P_2$ is any submodule of $P_1$ such that
the image in $P_1/P_2$ of $m$ is maximal.

\noindent (e) Let $M,N, \phi, m$ be as in (c) and assume $\phi(m)\ne 0$ . If $m$ is  a maximal
element of $M$,
 then $\phi(m)$ is a  maximal     element of $N$.

\medskip

\def\tzl{\tilde Z(\l)}
\def\tol{\tilde 1\ll}
\def\bzpln{\bold N_{\bold l}^n}
We shall denote by $\tzl$  the (baby) Verma module of $\tbu$ with
highest weight $\l$ and denote by $\tol$ a nonzero element in
$\tzl\ll$. Recall that to define $\Uxi$ we need to choose
$d_i\in\{1,2,3\}$ such that $(d_ia_{ij})$ is symmetric, where
$(a_{ij})$ is the concerned $n\times n$ Cartan matrix.
 Let $l_i$ be the order of $\xi^{2d_i}$. For $\l=(\l_1,...,\l_n)\in\bzn$
we set $\bold l\l=(l_1\l_1,...,l_n\l_n)$. We call $\l$ is $\bold
l$-restricted if $0\leq \l_i\leq l_i-1$ for all $i$. Denote by
$\bzpln$ the set of all $\bold l$-restricted elements in $\bzn$. The
following fact is well known.

\medskip

\noindent (f)  Let $\l\in\bzpln ,\ \l'\in\bzn $. Set $\mu=\l+\bold
l\l'\in\bzn $. Then $f_i^{(\l_i+1)}\tilde 1_\mu$ is maximal  in
$\tilde Z(\mu)$  if $\l_i\ne l_i-1$.

\bigskip

\section{ Some basic facts}

 \noindent{\bf 2.1.} From now on we assume that $\Uxi$ is of type
$B_2$. In this section we recall some basic facts about $\Uxi$ and
the Verma modules $\tzl$.  For completeness and fix notations,
  we give the definition
of $\Uxi$ and $\tzl$.

\def\bN{\bold N}
\def\bZ{\bold Z}
\def\ea{e_i^{(a)}}
\def\fa{f_i^{(a)}}
\def\kca{  \left[ \begin{array}{r }
      k_i,c \\ a\
             \end{array} \right] }
Let $a_{ii}=2$, $a_{12}=-2,\ a_{21}=-1$.  Let $U$ be the
associative algebra over $\bold Q(v)$ ($v$ an indeterminate) generated
by $e_i, f_i, k_i, $ $k_i^{-1}$ $ (i=1,2)$ with relations
$$k_1k_2=k_2k_1,\qquad k_ik_i^{-1}=k_i^{-1}k_i=1$$
$$k_ie_j=v^{ia_{ij}}e_jk_i,\qquad k_if_j=v^{-ia_{ij}}f_jk_i,$$
$$e_if_j-f_je_i=\delta_{ij}\frac{k_i-k_i^{-1}}{v_i-v_i^{-1}}$$
$$e_1e_2^2-(v^2+v^{-2})e_2e_1e_2+e_2^2e_1=0$$
$$e_1^3e_2-(v^2+1+v^{-2})e_1^2e_2e_1+(v^2+1+v^{-2})e_1e_2e_1^2-e_2e_1^3=0$$
 $$f_1f_2^2-(v^2+v^{-2})f_2f_1f_2+f_2^2f_1=0$$
$$f_1^3f_2-(v^2+1+v^{-2})f_1^2f_2f_1+(v^2+1+v^{-2})f_1f_2f_1^2-f_2f_1^3=0$$

\noindent where $v_1=v$ and $v_2=v^2$. Let $U'$ be the $A=\bold Z[v,
v^{-1}]$-subalgebra of $U$ generated by all $\ea=e_i^a/[a]_i!,
\fa=f_i^a/[a]_i!, k_i, k_i^{-1},\ a\in\bold N,\ i=1,2$, where
$[a]_i!=\prod_{h=1}^{a}\frac{v^{ih}-v^{-ih}}{v^i-v^{-i}}$ if
$a\geq 1$ and $[0]_i!=1$.
 Note that $\kca=\prod_{h=1}^{a}\frac{v_i^{c-h+1}\k-v_i^{-c+h-1}\kk}
 {v_i^h-v_i^{-h}}$ is in $U'$ for all
$c\in \bold Z,\ a\in\bold N$. We understand that $\kca=1$ if
$a=0$. Note that  $f'_{12}{}^{(a)}=(f_1f_2-v^2f_2f_1)^a/[a]!$ and
$f_{12}^{(a)}=(f_2f_1-v^2f_1f_2)^a/[a]!$ are in $U'$ for all
$a\in\bN$. Also
$$f'_{112}{}^{(a)}=\frac{(f_1f'_{12}-f'_{12}f_1)^a}{(v+v^{-1})^a[a]_2!}
\quad\text{and}\quad
f_{112}{}^{(a)}=\frac{(f_{12}f_1-f_1f_{12})^a}{(v+v^{-1})^a[a]_2!}$$
are in $U'$ for all $a\in\bN$. Regard $\qxi$ as an $A$-algebra by
specializing $v$ to $\xi$. Then $\Uxi=U'\otimes_A\qxi$. See [L3].

\def\e{e_i}
\def\f{f_i}
For convenience, the images in $\Uxi$ of $\ea,\fa,
f'_{12}{}^{(a)}, f_{12}^{(a)},f'_{112}{}^{(a)},
f_{112}^{(a)},\k,\kk,$ $\kca$ etc. will be denoted by the same
notation respectively. Let $l$ be the order of $\xi$ and $l_i$ be
the order of $\xi^{2i}$. In $\Uxi$ we have $\e^{l_i}=\f^{l_i}=0$.
{\bf For simplicity in this paper we assume that $l$ is odd. } Then
$l_1=l_2=l$.  The Frobenius kernel $\bu$ of $\Uxi$ is the
subalgebra of $\Uxi$ generated by all $\e,\f,\k,\kk,\ i=1,2$. Its
negative part $\bu^-$ is generated by all $\f$.   Note that
$f'_{12}{}^{(a)}, f_{12}^{(a)},f'_{112}{}^{(a)}, f_{112}^{(a)},$
are in $\bu^{-}$ if $ 0\leq a\leq l-1$. The subalgebra $\tbu$ of
$\Uxi$ is generated by all $\e,\f,\k,\kk,\kca,\ i=1,2;
c\in\bZ,a\in\bN$.

For $\l=(\l_1,\l_2)\in\bZ^2$, we denote by $\tilde I_\l$
 the left ideal of $U_\xi$ generated by all
$\e^{(a)}\ (a>0), \k-\xi^{i\l_i},\kca-{\left[\begin{array}{r}
\l_i+c\\ a\
\
\
\end{array}\right] }_{\xi^i}$. (We denote by
${\left[\begin{array}{r} b\\ a\end{array}\right]}_{\xi^i}$ the value
at $\xi^i$ of $\prod_{h=1}^{a}\frac{v^{b-h+1}-v^{
b+h-1}}{v^h-v^{-h}}$ for any $b\in\bZ$ and $a\in\bN$.) The Verma
module $Z(\l)$ of $\Uxi$ is defined to be $U_\xi/\tilde I_\l$. Let
$\tilde 1\ll$ be the image in $Z(\l)$ of 1.  The Verma module $\tzl$
of $\tbu$ is defined to be the $\tbu$-submodule of $Z(\l)$ generated
by $\tilde 1\ll$ . Given non-negative integers $a$ and $b$, we set
$$x_{a,b}=f_1^{(a)}f_2^{(a+b)}f_1^{(a+2b)}f_2^{(b)}=
f_2^{(b)}f_1^{(a+2b)}f_2^{(a+b)}f_1^{(a)}.$$ Recall that  $l$ is
the order of $\xi$. The following result is a special case of [X1,
4.2 (ii)].

\def\bzpln{\bold N_{\bold l}^2}
\def\bzn{\bZ^2}
\noindent(a) Assume $0\leq a,b\leq l-1$, $c,d\in\bZ$, and let
$\mu=(lc-1+a , ld-1+b )$. Then the element $x_{a,b}$ is in $\bu^-$
and $x_{a,b} \tilde 1_\mu$ is maximal in $\tilde Z(\mu)$ and
generates the unique irreducible submodule of $\tilde Z(\mu)$. The
irreducible submodule is isomorphic to  $\tilde L(lc-1-a ,ld-1-b
)$.

The argument for [X1, 4.4(iv)]also gives the following result.

\noindent(b) Keep the assumption and notations in (a). Let
$p,q,s,t\in \bold N$ such that
$x=f_1^{(a+2b-pl)}f_2^{(a+b-ql)}f_1^{(a)}$ and
$y=f_2^{(a+b-sl)}f_1^{(a+2b-tl)}f_2^{(b)}$ are nonzero elements,
then $E_ix\tilde 1_\mu=E_ix\tilde 1_\mu=0$ for $i=1,2$. If $x$ and
$y$ are further in $\bu^-$, then $x\tilde 1_\mu$ and $y\tilde 1_\mu$
are maximal in $\tilde Z(\mu)$.

We shall need a few formulas, which are due to Lusztig (see [L3, L4]. In $\Uxi$
we have

\medskip

\def\no{\noindent}
\noindent(c)
$\displaystyle{f_i^{(a)}\f^{(b)}=\left[\begin{array}{r}a+ b\\ a\ \
\end{array}\right]_{\xi^i}\f^{(a+ b)}},$

\no (d) $f^{(i)}_{12} f^{(j)}_2 = \xi^{2ij} f^{(j)}_2 f^{(i)}_{12},
$
\medskip
\def\N{\bold N}
\def\st{\stackrel}
\def\sc{\scriptstyle}
\newcommand\ds{\displaystyle\sum}

\no (e) $f^{(i)}_{112} f^{(j)}_{12} = \xi^{2ij} f^{(j)}_{12}
f^{(i)}_{112}, $ \medskip

\no (f) $ f^{(i)}_1 f^{(j)}_{112} = \xi^{2ij} f^{(j)}_{112}
f^{(i)}_1, $ \medskip

\no (g) $ f^{(i)}_2f^{(j)}_{112}  = \ds _{ \st{\sc r, s, t \in \N}
{\st{\sc r +s = j}  {s+t = i } } } \xi ^{-2 r s- 2 s
t}\displaystyle\prod ^ s_{h=1} (\xi^{-4h+2} -1)f^{(r)}_{112}
f^{(2s)} _{12}f^{(t)}_2 , $
\medskip

\no (h) $ f^{(i)}_{12}f^{(j)}_1  = \ds _{ \st{\sc r,s,t \in \N} {
\st{\sc r+s= j} {s+t = i }}  } \xi^{-rs- st+ s}\displaystyle\prod ^
s _{h=1} (\xi^{-2h}+1)  f^{(r)}_1 f^{(s)}_{112}f^{(t)}_{12}, $
\medskip

\no (i) $ f^{(i)}_2 f^{(j)}_1 = \ds _{ \st{\sc r, s, t, u \in\N}
{\st{\sc s+t+u=i} {r+2s+t = j }} } \xi^{2 r u + 2 su+ rt} f^{(r) }_1
f^{(s)}_{112} f^{(t)}_{12} f^{(u)}_2.$

\medskip

\no (j) $ f'{}^{(i)}_{112} f^{(j)}_2 = \ds _{ \st{\sc r, s, t \in
\N} { \st{\sc r +s = j} {s+t = i}  } } \xi ^{-2 r s- 2 s
t}\displaystyle\prod ^ s_{h=1} (\xi^{-4h +2} -1) f^{(r)}_2
f'{}^{(2s)} _{12} f'{}^{(t)}_{112}, $ \medskip

\no (k) $ f^{(i)}_1 f'{}^{(j)}_{12} = \ds _{ \st{\sc r,s,t \in \N}
{\st{\sc r+s= j} {s+t = i}  } } \xi^{-rs- st+ s}\displaystyle\prod ^
s _{h=1} \xi^{-2h}+1) f'{}^{(r)}_{12} f'{}^{(s)}_{112} f^{(t)}_1, $
\medskip

\no (l) $ f^{(i)}_1 f^{(j)}_2 = \ds _{ \st{\sc r, s, t, u \in\N}
{\st{\sc r+ s+t= j} {s+ 2 t + u = i } } } \xi^{2 r u + 2 r t+us}
f^{(r) }_2 f'{}^{(s)}_{12} f'{}^{(t)}_{112} f^{(u)}_1.$

\def\lan{\langle}
\def\ran{\rangle}

\noindent(m) Assume $0\leq a_0,b_0\leq l-1$ and $a_1,b_1\in\bZ$. We
have $${\left[\begin{array}{r} {a_0+a _1l}\\ {b_0+b
_1l}\end{array}\right]}_{\xi^i}=   {\left[\begin{array}{r} {a_0}\\
{b_0}\end{array}\right]}_{\xi^i}{\left(\begin{array}{r} {a_1}\\
{b_1}\end{array}\right)},$$ where ${\left(\begin{array}{r} {a_1}\\
{b_1}\end{array}\right)}$ is the ordinary binomial coefficient.

Using (i), (l) and (m) we get

\noindent(n) If $0\le a\le l-1$, then
$f_1^{(l)}f_2^{(a)}-f_2^{(a)}f_1^{(l)}$ and $
f_2^{(l)}f_1^{(a)}-f_1^{(a)}f_2^{(l)}$ are in $\bu^-$.

\noindent(o) Let $0\le a,b,c\le l-1$. Then
$f_1^{(a)}f_2^{(b)}f_1^{(c)}=0$ and $f_1^{(a)}f_2^{(l+b)}f_1^{(c)}$
is in $\bu^-$, if $a+c-2b\ge l$. Similarly
$f_2^{(a)}f_1^{(b)}f_2^{(c)}=0$ and $f_2^{(a)}f_1^{(l+b)}f_2^{(c)}$
is in $\bu^-$,
  if $a+c-b\ge l$

  The assertions (n) and (o) will be frequently used in
  computations.

Let $\a_1=(2,-1),\a_2=(-2,2)\in\bZ^2$. The set of positive roots
is $R^+=\{\a_1,\a_2,\a_1+\a_2,2\a_1+\a_2\}$. Let $W$ be the Weyl
group generated by the simple reflections $s_i$ corresponding to
$\a_i$. Assume that $l\ge 5$. Then $\lan \rho,\b^\vee\ran< l$ for
all $\b\in R^+$, where $\rho=(1,1)$. For $\lambda,\mu\in\bZ^2$, we write
that $\lambda\leq\mu$ if $\mu-\lambda=a\a_1+b\a_2$ for some non-negative
integers $a,b$.

We say that $x\in \bu^-$ is homogeneous (of degree $\beta$) if there
exists $\beta\in\bZ^2$ such that $\kca x=x{\left[\begin{array}{r}
k_i,c+\lan\beta,\a_i^\vee\ran\\ a\ \ \
\end{array}\right]} $ and $k_ix=\xi^{i\lan\beta,\a_i^\vee\ran}xk_i$ for all $c\in\bZ$ and $a\in\bN.$

\bigskip

\def\tol{\tilde 1_\l}
\def\t1m{\tilde 1_\mu}

  \noindent {\bf 2.2.} The $W$-orbit of $\l=(0,0)$ (dot
action)  consists of the following 8 elements,

$\l,\ \  s_1.\l=\l- \a_1,\ \  s_2.\l=\l- \a_2,\ \
 s_2s_1.\l=\l- \a_1- 2\a_2,\ \ s_1s_2.\l=\l- 3\a_1- \a_2,\ \
s_1 s_2s_1.\l=\l-4\a_1- 2\a_2,$ $s_2s_1s_2.\l=\l- 3\a_1- 3\a_2, \
\ s_1s_2s_1s_2.\l=\l-4\a_1- 3\a_2.$

Let $a,b$ be integers and $\l=(la,lb)$. Using 1.1  (e-f) and 2.1
(a-b) we get

\noindent(1)  The following  elements are maximal     in $\tilde
Z(\l)$: $$\tol,\qquad f_1 \tol,\qquad f_2 \tol,\qquad f_2^{(2)}f_1
\tol,\qquad
 f_1^{(3)}f_2 \tol,
 $$ $$
 f_2^{(2)}f_1^{(3)}f_2 \tol,\qquad
  f_1^{(3)}f_2^{(2)}f_1 \tol,\qquad
  f_1  f_2^{(2)}f_1^{(3)}f_2 \tol.$$

\noindent(2) Let $\mu=\l+(l -1)\a_1$. The following elements are
maximal     in $\tilde Z(\mu)$: $$\t1m,\qquad f_1^{(l
-1)}\t1m,\qquad f_2^{(2)}\t1m,\qquad f_2 f_1^{(l -1)}\t1m,\qquad
f_1^{(3)}f_2^{(2)}\t1m,$$ $$ f_2 f_1^{(3)}f_2^{(2)}\t1m, \qquad
f_1^{(3)}f_2^{(l +1)}f_1^{(l -1)}\t1m,$$ $$ f_2^{(2)}f_1^{(l
+3)}f_2^{(l +1)}f_1^{(l -1)}\t1m.$$

 \noindent(3)
 Let $\mu=\l+(l -1)\a_2$. The
following elements are maximal in $\tilde Z(\mu)$: $$\t1m,\qquad
f_1^{(3)}\t1m,\qquad f_2^{(l -1)}\t1m,\qquad f_2^{(2)}f_1^{(
3)}\t1m,\qquad f_1 f_2^{(l -1)}\t1m,$$ $$ f_1 f_2^{(2)}f_1^{(
3)}\t1m,\qquad f_2^{(2)}f_1^{(l +1)}f_2^{(l -1)}\t1m,\qquad f_1^{
(3)}f_2^{(l +2)}f_1^{(2l +1)}f_2^{(l -1)}\t1m.$$

\noindent(4) Let $\mu=\l+(l -1)\a_1+(l  -2)\a_2$. The following
elements are maximal     in $\tilde Z(\mu)$: $$\t1m,\qquad
f_1^{(3)}\t1m,\qquad f_2^{(l  -2)}\t1m,\qquad f_2 f_1^{(3)}\t1m,$$
$$ f_1^{(l -1)}f_2^{(l  -2)}\t1m,\qquad
 f_2 f_1^{(l -1)}f_2^{(l  -2)}\t1m.$$
 $$f_1^{(l -1)} f_2^{(l +1)} f_1^{(3)}\t1m,\qquad
 f_1^{(3)}f_2^{(l +1)}f_1^{(2l -1)}f_2^{(l  -2)}\t1m.$$

\noindent(5) Let $\mu=\l+(l  -3)\a_1+(l -1)\a_2$. The following
elements are maximal     in $\tilde Z(\mu)$: $$\t1m,\qquad f_1^{(l
 -3)}\t1m,\qquad f_2^{(2)}\t1m,\qquad f_2^{(l -1)}f_1^{(l
 -3)}\t1m,$$ $$ f_1 f_2^{(2)}\t1m,\qquad
 f_1 f_2^{(l -1)}f_1^{(l  -3)}\t1m,$$
 $$
f_2^{(l -1)}f_1^{(l +1)}f_2^{(2)}\t1m,\qquad f_2^{(2)} f_1^{(l
+1)}f_2^{(l -1)}f_1^{(l  -3)}\t1m.$$

\noindent(6) Let $\mu=\l+(l- 4)\a_1+(l  -2)\a_2$. The following
elements are maximal     in $\tilde Z(\mu)$: $$\t1m,\qquad f_1^{(l
 -3)}\t1m,\qquad f_2 \t1m,\qquad f_2^{(l
 -2)}f_1^{(l  -3)}\t1m,$$ $$ f_1^{(l -1)}f_2 \t1m,\qquad
 f_1^{(l -1)}f_2^{(l  -2)}f_1^{(l  -3)}\t1m.$$
  $$
f_2^{(l  -2)}f_1^{(l -1)}f_2 \t1m,\qquad f_2 f_1^{(l -1)}f_2^{(l
-2)}f_1^{(l  -3)}\t1m.$$

\noindent (7) Let $\mu=\l+(l  -3)\a_1+(l  -3)\a_2$. The following
elements are maximal     in $\tilde Z(\mu)$: $$\t1m,\qquad f_1
\t1m,\qquad f_2^{(l  -2)}\t1m,\qquad f_2^{(l -1)}f_1 \t1m,$$ $$
f_1^{(l  -3)}f_2^{(l
 -2)}\t1m,\qquad
 f_1^{(l  -3)}f_2^{(l -1)}f_1 \t1m.$$
  $$
f_2^{(l -1)}f_1^{(2l  -3)}f_2^{(l  -2)}\t1m,\qquad f_2^{(l
 -2)} f_1^{(2l  -3)}f_2^{(l -1)}f_1 \t1m.$$

\no(8) Let $\mu=\l+(l  -4)\a_1+(l  -3)\a_2$. The following
elements are maximal     in $\tilde Z(\mu)$: $$\t1m,\qquad f_1^{(l
-1)}\t1m,\qquad f_2^{(l -1)}\t1m,\qquad f_2^{(l  -2)}f_1^{(l
-1)}\t1m,$$ $$ f_1^{(l  -3)}f_2^{(l -1)}\t1m,\qquad
 f_1^{(l  -3)}f_2^{(l  -2)}f_1^{(l -1)}\t1m.$$
  $$
f_2^{(l  -2)}f_1^{(2l  -3)}f_2^{(l -1)}\t1m,\qquad f_2^{(l -1)}
f_1^{(3l  -3)}f_2^{(2l  -2)}f_1^{(l-1 )}\t1m.$$

\section{maximal and primitive elements  of $\tzl$ for type $B_2$}

In this section we determine the maximal and primitive elements in
$\tzl$ (or equivalently in any
 highest weight module of $\tbu$).
  To avoid complicated expressions we only work
 with some weights in the $W$-orbit of $(0,0)$.
 For general cases the approach is completely similar. Throughout
 the paper $l$  is odd and is greater than or equal to 5.

\bigskip

\noindent {\bf Theorem 3.1.} {\sl   Let $a,b$ be integers and
$\l=(la,lb) $. Then}

\noindent(i){\sl  The following 8 elements are maximal     in
$\tilde Z(\l)$: $$\tol,\qquad f_1 \tol,\qquad f_2 \tol,\qquad
f_2^{(2)}f_1 \tol,\qquad
 f_1^{(3)}f_2 \tol,
 $$ $$
 f_2^{(2)}f_1^{(3)}f_2 \tol,\qquad
  f_1^{(3)}f_2^{(2)}f_1 \tol,\qquad
  f_1  f_2^{(2)}f_1^{(3)}f_2 \tol.$$}
\noindent(ii){\sl  The following 12 elements are  primitive
elements in $\tilde Z(\l)$ but not maximal:
$$[f_1^{(3)},f_2^{(l)}]f_2\tol,\quad\frac{x_{l-1,2}}{f_1^ {(l-1)}
} \tol,\quad[f_2^{(2)},f_1^{(l)}]f_1^{(3)}f_2 \tol$$
$$[f_2^{(2)},f_1^{(l)}]f_1\tol,\quad\frac{x_{3,l-1}}{f_2^{(l-1)} }
\tol,\quad\frac{x_{3,l-1}}{f_1^{(l)}f_2^{(l-1)} }
\tol,\quad[f_1^{(3)},f_2^{(l)}]f_2^{(2)}f_1 \tol,$$
$$\frac{x_{3,l-2}}{f_1^{(l-1)}f_2^{(l-2)} } \tol,\quad
f_2^{(l-1)}f_1^{(l)}f_2\tol,\quad
f_1f_2^{(l-1)}f_1^{(l)}f_2\tol,$$
$$\frac{x_{l-1,l-1}}{f_1^{(l-3)}f_2^{(l-2)}f_1^{(l-1)}}f_2\tol,\quad
\frac{x_{l-1,l-1}}{f_2^{(l-2)}f_1^{(2l-3)}f_2^{(l-1)}}f_1\tol.$$
(See 1.1 for the definition of $x_{a,b}$.
 Convention: $[x,y]=xy-yx$ and $\frac xy$
 stands for an arbitrary homogeneous element $z$ in $\bu^-$ such that $zy=x$.)

Moreover  no maximal     element in $\tilde Z(\l)$  has the same
weight as any of the above 12 elements.}

\def\tl{\tilde L}

\noindent(iii) {\sl The maximal and primitive          elements in
(i-ii) provide 20 composition factors of $\tilde Z(\l)$, which are
$$\tilde L(\l),\qquad\tilde L(\l- \a_1),\qquad\tilde L(\l-
\a_2),$$
$$\tilde L(\l- \a_1- 2\a_2),\qquad\tilde L(\l- 3\a_1-
\a_2)$$
 $$\tilde L(\l- 3\a_1- 3\a_2),\qquad\tilde L(\l- 4\a_1-
2\a_2),\qquad\tilde L(\l- 4\a_1- 3\a_2),$$ $$
\tl(\l-3\a_1-(l+1)\a_2),\qquad\tl(\l-(l+3)\a_1-(l+3)a_2),\qquad\tl(\l-(l+3)\a_1-3\a_2),
$$
$$\tl(\l-(l+1)\a_1-2\a_2),\qquad\tl(\l-(2l+4)\a_1-(l+2)\a_2),\qquad\tl(\l-(l+4)\a_1
-(l+2)\a_2),$$ $$
 \tl(\l-4\a_1-(l+2)\a_2),\qquad\tl(\l-(l+3)\a_1-(l+1)\a_2),\qquad\tl(\l-l
 \a_1-l\a_2),$$ $$
 \tl(\l-(l+1)\a_1-l\a_2),\qquad\tl(\l-2l\a_1-2l\a_2),
 \qquad\tl(\l-2l\a_1-l\a_2).$$
 Moreover, $\tilde Z(\l)$ has
only the 20 composition factors.}

{
\def\g{\gamma}

\noindent{\it Proof:} (i) According to 1.1 (e-f), we see that (i)
is true.

\def\g{\gamma}

\noindent(ii) Now we argue for (ii).
\def\vp{\varphi}

 \noindent(1)
Consider the homomorphism: $$\varphi_1:\  \tilde Z(\l)\to\tilde
Z(\l+(l -1)\a_1)),\quad \tol\to t_1=f_1^{(l -1)}\tilde 1_{\l+(l
-1)\a_1}.$$

Let $$x_1=(f_1^{(3)}f_2^{(l)}- f_2^{(l)}f_1^{(3)})f_2 \in\uxin.$$
Note that $$f_1^{(3)}f_2^{(l +1)}t_1 =f_1^{(3)}f_2^{(l +1)}f_1^{(l
-1)}\tilde 1_{\l+(l -1)\a_1}= x_1t_1.$$ Using 1.1 (c) we see that
$x_1\tol$ is a primitive element of weight $\g_{38}=\l- 3\a_1-(l
+1)\a_2$.

  Let
$$y_1=(f_2^{(2)}f_1^{(l)}- f_1^{(l)}f_2^{(2)})x_1\in\uxin.$$ Note
that $f_2^{(2)}f_1^{(3)}f_2^{(l +1)}f_1^{(l -1)}=0$. We then can
check that $$f_2^{(2)}f_1^{(l+3)}f_2^{(l +1)}t_1 =f_2^{(2)}f_1^{(l
+3)}f_2^{(l +1)}f_1^{(l -1)}\tilde 1_{\l+(l -1)\a_1}=y_1t_1.$$ Using
1.1 (c) we see that $y_1\tol$ is a  primitive element of weight
$\g_{48}=\l-(l +3)\a_1-(l +3)\a_2$. Note that we have
$y_1=\frac{x_{l-1,2}}{f_1^{(l-1)}}.$

\def\tbu{\tilde{\bold u}}
Note that $f_2^{(l +1)}f_1^{(l -1)}=f_2f_2^{(l)}f_1^{(l -1)}$, so we
have $xf_2^{(l +1)}f_1^{(l -1)}=0$ if $xf_2=0$ and $x\in\bu^-$. Thus
we have a homomorphism (recall that $\tilde 1_{\l+(l -1)\a_1+l\a_2}$
is also an element of the Verma module $Z({\l+(l -1)\a_1+l\a_2})$ of
$U_\xi$):
$$\psi_1:\ \tbu f_2\tol\to\tbu f_2^{(l +1)}f_1^{(l
-1)}\tilde 1_{\l+(l -1)\a_1+l\a_2},$$
$$\quad f_2 \tol\to f_2^{(l +1)}f_1^{(l -1)}\tilde 1_{\l+(l -1)\a_1+l\a_2}.$$

Note that $\psi_1(f_1^{(3)}f_2 \tol)= f_1^{(
 3)}f_2^{(l +1)}f_1^{(l -1)}\tilde 1_{\l+(l
-1)\a_1+l\a_2}$ and $f_1^{(
 3)}f_2^{(l +1)}f_1^{(l -1)}$ is in $\bu^-$.

Let $z_1=f_2^{(2)}f_1^{(l)}- f_1^{(l)}f_2^{(2)}\in\uxin$. Using 1.1
(c) we see that $z_1f_1^{(3)}f_2 \tol$ is a primitive element of
weight $\g_{47}=\l-(l +3)\a_1- 3\a_2$. Note that
$z_1f_1^{(3)}f_2=\frac{x_{l-1,2}}{f_2^{(l)}f_1^{(l-1)}}.$

\medskip

 \noindent(2)
Now we consider the homomorphism: $$\vp_2:\ \tilde Z(\l)\to\tilde
Z(\l+(l -1)\a_2)),\quad \tol\to t_2=f_2^{(l -1)}\tilde 1_{\l+(l
-1)\a_2}.$$

Let $$x_2=(f_2^{(2)}f_1^{(l)}-f_1^{(l)}f_2^{(2)})f_1\in\uxin.$$ Note
that $$f_2^{(2)}f_1^{(l +1)}t_2 =f_2^{(2)}f_1^{(l +1)}f_2^{(l
-1)}\tilde 1_{\l+(l -1)\a_2}=x_2t_2.$$ Using 1.1 (c) we see that
$x_2\tol$ is a  primitive element of weight $\g_{31}=\l-(l +1)\a_1-
2\a_2$.

 Let $y_2$ be homogeneous in $\uxin$ such that
 $y_2t_2=f_1^{(3)}f_2^{(l +2)}f_1^{(2l +1)}f_2^{(l -1)}\t1m$, here $\mu=\l+(l-1)\a_2$.
 According to 1.1 (c) we know that $y_2\tol$ is a primitive
 element of weight $\g_{41}=\l-(2l+4)\a_1-(l +2)\a_2$. Note that
$y_2=\frac{x_{3,l-1}}{ f_2^{(l-1)}}.$

As the reason for $\psi_1$, we have a homomorphism: $$\psi_2:\ \tbu
f_1\tol\to\tbu f_1^{(l +1)}f_2^{(l -1)}\tilde 1_{\l+l\a_1+(l
-1)\a_2},$$
$$\quad f_1 \tol\to  f_1^{(l +1)}f_2^{(l
-1)}\tilde 1_{\l+l\a_1+(l -1)\a_2}.$$

Note that $\psi_2(f_2^{(2)}f_1 \tol)= t'_2= f_2^{(2)}f_1^{(l
+1)}f_2^{(l -1)}\tilde 1_{\l+l\a_1+(l -1)\a_2}$ and
$f_2^{(2)}f_1^{(l +1)}f_2^{(l -1)}$ is in $\bu^-$. Let $z_2$ be
homogeneous in $\uxin$ such that $$z_2t'_2
=f_1^{(3)}f_2^{(l+2)}f_1^{(2l +1)}f_2^{(l -1)}\tilde 1_{\l+l\a_1+(l
-1)\a_2}. $$ Using 1.1 (c) we see that $z_2f_2^{(2)}f_1 \tol$ is a
primitive element of weight $\g_{42}=\l-(l+ 4)\a_1-(l +2)\a_2$. Note
that $z_2f_2^{(2)}f_1=\frac{x_{3,l-1}}{f_1^{(l)}f_2^{(l-1)}}.$

 We also have a homomorphism (recall that $\tilde 1_{\l+2l\a_1+(l
-1)\a_2}$ is also an element of the Verma module $Z({\l+2l\a_1+(l
-1)\a_2})$ of $U_\xi$) : $$\theta_2:\ \tbu f_1\tol\to\tbu f_1^{(2l
+1)}f_2^{(l -1)}\tilde 1_{\l+2l\a_1+(l -1)\a_2},$$
$$\quad  f_1 \tol\to t''_2= f_1^{(2l +1)}f_2^{(l
-1)}\tilde 1_{\l+2l\a_1+(l -1)\a_2}.$$

Let $w_2=[f_1^{(3)},f_2^{(l)}]\in\uxin$. Then $$w_2f_2^{(2)}t''_2
=f_1^{(3)}f_2^{(l+2)}f_1^{(2l +1)}f_2^{(l -1)}\tilde 1_{\l+2l\a_1+(l
-1)\a_2}. $$ Using 1.1 (c) we see that $w_2f_2^{(2)}f_1 \tol$ is a
primitive element of weight $\g_{44}=\l-4\a_1-(l +2)\a_2$. Note that
$w_2f_2^{(2)}f_1 =\frac{x_{3,l-1}}{f_1^{(2l)}f_2^{(l-1)}}.$

\medskip

\noindent(3) Now we consider the homomorphism: $$\tilde
Z(\l)\to\tilde Z(\l+(l -1)\a_1+(l -2)\a_2),$$ $$ \tol\to
t_3=f_1^{(l -1)}f_2^{(l
 -2)}\tilde 1_{\l+(l -1)\a_1+(l -2)\a_2}.$$

Let $x_3$ be homogeneous in $\uxin$ such that
$$x_3t_3=f_1^{(3)}f_2^{(l +1)}f_1^{(2l -1)}f_2^{(l  -2)}\tilde
1_{\l+(l -1)\a_1+(l
 -2)\a_2}.$$ Using 1.1 (c) we know that $x_3\tol$ is a primitive
element of weight $\g_{37}=\l-(l +3)\a_1-(l +1)\a_2$.  Note that
$x_3=\frac{x_{3,l-2}}{f_1^{(l-1)}f_2^{(l-2)}}.$

\noindent(4) Since $ f_1^{(l  -3)}f_2^{(l  -2)}\in f_2 \uxin$ (see
2.1), we have a surjective homomorphism $$\varphi_4:\ \tbu f_2
\tol\to \tbu f_1^{(l
 -3)}f_2^{(l  -2)}\tilde 1_{\l+(l  -3)\a_1+(l  -3)\a_2}$$ $$ f_2
\tol\to
 f_1^{(l  -3)}f_2^{(l  -2)}\tilde
1_{\l+(l  -3)\a_1+(l  -3)\a_2}.$$

Let $x_4=f_2^{(l -1)}f_1^{(l)}- f_1^{(l)}f_2^{(l -1)}\in \uxin.$
Then $$x_4f_2 \tol=f_2^{(l -1)}f_1^{(l)}f_2 \tol$$ is a primitive
element of weight $\g_{34}=\l-l\a_1-l\a_2$.

Let $y_4=f_1 x_4f_2$. Then $y_4\tol$ is a primitive element of
                          weight $\g_{45}=\l-(l +1)\a_1-l\a_2$.

\medskip

\no(5) Consider the homomorphism $$\tbu f_2 \tol\to \tbu f_1^{(l
 -3)}f_2^{(l  -2)}f_1^{(l -1)}\tilde 1_{\l+(2l  -4)\a_1+(l
 -3)\a_2}$$ $$ f_2 \tol\to
 f_1^{(l  -3)}f_2^{(l  -2)}f_1^{(l -1)}\tilde
1_{\l+(2l  -4)\a_1+(l  -3)\a_2}.$$ Let $x_5$ be homogeneous in
$\uxin$ such that $$x_5f_1^{(l  -3)}f_2^{(l  -2)}f_1^{(l
-1)}=f_2^{(l -1)} f_1^{(3l -3)}f_2^{(2l  -2)}f_1^{(l -1)}.$$ By 1.1
(c), $x_5f_2\tol$ is primitive and is of weight
$\g_{35}=\l-2l\a_1-2l\a_2$. Note that
$x_5=\frac{x_{l-1,l-1}}{f_1^{(l-3)}f_2^{(l-2)}f_1^{(l-1)}}.$

Consider the surjective homomorphism  $$\tbu f_1 \tol\to \tbu
f_2^{(l  -2)}f_1^{(2l  -3)}f_2^{(l -1)}\tilde 1_{\l+(2l
 -4)\a_1+(2l  -3)\a_2}$$ $$  f_1 \tol\to
  f_2^{(l  -2)}f_1^{(2l  -3)}f_2^{(l -1)}\tilde
1_{\l+(2l  -4)\a_1+(2l  -3)\a_2}.$$ Let $y_5$ be homogeneous in
$\uxin$ such that $$y_5 f_2^{(l  -2)}f_1^{(2l  -3)}f_2^{(l
-1)}=f_2^{(l -1)} f_1^{(3l  -3)}f_2^{(2l
 -2)}f_1^{(l-1  )}.$$  By 1.1 (c), $y_5f_1\tol$ is primitive and is of
weight $\g_{32}=\l-2l\a_1-l\a_2$. Note that
$y_5=\frac{x_{l-1,l-1}}{f_2^{(l-2)}f_1^{(2l-3)}f_2^{(l-1)}}.$

We may also consider the homomorphism: $$\tbu f_2 \tol\to\tbu
f_1^{(l-3)}f_2^{(2l-2)}f_1^{(l-1)}
 \tilde
1_{\l+(2l-4)\a_1+(2l-2)\a_2} ,$$ $$  \tbu f_1 \tol\to
f_1^{(l-3)}f_2^{(2l-2)}f_1^{(l-1)}
 \tilde
1_{\l+(2l-4)\a_1+(2l-2)\a_2} .$$   Using 1.1 (c) we know that
$[[f_2^{(l-1)},f_1^{(l)}],f_1^{(l)}]f_2\tol=\frac{x_{l-1,l-1}}{f_1^{(l-3)}f_2^{(2l-2)}f_1^{(l-1)}}f_2\tol$
is also a primitive element of weight $\l-2l\a_1-l\a_2$.

\bigskip

The element $f_1f_2^{(2)}f_1^{(3)}f_2\tol$ generates the unique
irreducible submodule of $\tzl$. Clearly, the weight of any element
in (ii) is not greater than  $\l-4\a_1-3\a_2$, therefore  no maximal
element in $\tzl$ has the same weight as any of the elements in
(ii).

\noindent(iii) Using (i), (ii) and  1.1 (b), we see that $\tzl$ has
the 20 composition factors. The dimensions of irreducible $\tilde
{\bold u}$-modules are known (see [APW]). By a comparison of
dimensions we know that $\tzl$ has only the 20 composition factors.

The theorem is proved.

\bigskip

\noindent {\bf Theorem 3.2.} {\sl   Let $a,b$ be integers and
$\l=(la ,lb-3) $. Then}

\noindent(i){\sl  The following  elements are maximal     in
$\tilde Z(\l)$: $$\tol,\qquad f_1 \tol,\qquad f_2^{(l
-2)}\tol,\qquad f_2^{(l -1)}f_1 \tol,$$ $$ f_1^{(l  -3)}f_2^{(l
 -2)}\tol,\qquad
 f_1^{(l  -3)}f_2^{(l -1)}f_1 \tol.$$
  $$
f_2^{(l -1)}f_1^{(2l-3)}f_2^{(l-2)}\tol,\qquad f_2^{(l
 -2)} f_1^{(2l-3)}f_2^{(l -1)}f_1 \tol,$$}
$$\frac{f_2^{(l-1)}f_1^{(l-3)}}{f_2^{(2)}}\tol,\qquad
\frac{f_1f_2^{(l-1)}f_1^{(l-3)}}{f_2^{(2)}}\tol,\quad
f_1^{(l-3)}[f_2^{(l-1)},f_1^{(l)}]f_1\tol.$$

\noindent(ii){\sl  The following   elements are  primitive
elements in $\tilde Z(\l)$ but not maximal:
$$\frac{x_{l-1,l-1}}{f_1^{(l-1)}}\tol,\quad\frac{x_{l-1,l-1}}{f_2^{(l)}
f_1^{(l-1)}}\tol,\qquad
[f_2^{(l-1)},f_1^{(l)}]f_1\tol,$$ $$
\frac{x_{3,l-1}}{f_2^{(2)}f_1^{(3)}}\tol, \quad\quad \frac{x_{3,l-1}}{f_2^{(l+2)}
f_1^{(3)}}\tol,\quad\quad f_1^{(l-1)}f_2^{(l)}f_1\tol,$$ $$
\frac{x_{3,l-2}}{f_2f_1^{(3)}}f_1\tol,
\quad[f_2^{(l-2)},f_1^{(l)}]\tol,\quad\frac{x_{l-1,2}}{f_2
f_1^{(3)}f_2^{(2)}}f_1\tol.$$

Moreover  no maximal     element in $\tilde Z(\l)$ has the same
weight as any of the above 9 elements.}

\def\tl{\tilde L}

\noindent(iii) {\sl The maximal and primitive          elements in
(i-ii) provide 20 composition factors of $\tilde Z(\l)$, which are
$$\tilde L(\l),\qquad\tilde L(\l- \a_1),\qquad\tilde L(\l-(l-2)
\a_2),$$ $$\tilde L(\l- \a_1-(l-1)\a_2),\qquad\tilde
L(\l-(l-3)\a_1-(l-2) \a_2)$$
 $$\tilde L(\l- (l-2)\a_1-(l-1)\a_2),\qquad\tilde L(\l-(2l-3)\a_1-
(2l-3)\a_2),\qquad\tilde L(\l-(2l-2)\a_1-(2l-3)\a_2),$$ $$
\tl(\l-(l-3)\a_1-(l-3)\a_2),\qquad\tl(\l-(l-2)\a_1-(l-3)a_2),
\qquad\tl(\l-(3l-3)\a_1-(3l-3)\a_2), $$
$$\tl(\l-(3l-3)\a_1-(2l-3)\a_2),\qquad\tl(\l-(l+1)\a_1-(l-1)\a_2),
\qquad\tl(\l-(2l-2)\a_1 -(l-1)\a_2),$$ $$
 \tl(\l-(2l+1)\a_1-(2l-1)\a_2),\qquad\tl(\l-(2l+1)\a_1-(l-1)\a_2),
 \qquad\tl(\l-l
 \a_1-l\a_2),$$ $$
 \tl(\l-l\a_1-(2l-2)\a_2),\qquad\tl(\l-l\a_1-(l-2)\a_2),
 \qquad\tl(\l-2l\a_1-l\a_2).$$
 Moreover, $\tilde Z(\l)$ has
only the 20 composition factors.}

\def\g{\gamma}

\noindent{\it Proof:} (i) According to 1.1 (e-f), we see that the first
8 elements in (i) are maximal.

Consider the homomorphism: $$\tilde Z(\l)\to\tilde
Z(\l+2\a_2),\quad \tol\to t_1=f_2^{(2)}\tilde 1_{\l+2\a_2}.$$

Since $f_2^{(l-1)}f_1^{(l-3)}$ is in $\bu^-f_2^{(2)}$, using 1.1 (c)
we see that $\frac{f_2^{(l-1)}f_1^{(l-3)}}{f_2^{(2)}}\tol$ and $
\frac{f_1f_2^{(l-1)}f_1^{(l-3)}}{f_2^{(2)}}\tol$ are primitive
elements of weights $\l-(l-3)\a_1-(l-3)\a_2$ and
$\l-(l-2)\a_1-(l-3)\a_2$ respectively. One can check directly that
the two elements are maximal. We will show that the last  element in
(i) is maximal in part (2) of the argument for (ii).

\noindent(ii) Now we argue for (ii).

 \noindent(1)
Consider the homomorphism: $$\tilde Z(\l)\to\tilde Z(\l+(l
-1)\a_1)),\quad \tol\to t_1=f_1^{(l -1)}\tilde 1_{\l+(l
-1)\a_1}.$$ Let $x_1=\frac{x_{l-1,l-1}}{f_1^{(l-1)}}\in\uxin$.
Using 1.1 (c) we see that $x_1\tol$ is a primitive element of
weight $\l-(3l-3)\a_1-(3l -3)\a_2$.

Consider the homomorphism: $$\tbu f_2^{(l-2)}\tol\to\tbu
f_2^{(2l-2)}f_1^{(l-1)}\tilde 1_{\l+(l -1)\a_1+l\a_2},$$ $$
f_2^{(l-2)}\tol\to  f_2^{(2l-2)}f_1^{(l-1)}\tilde 1_{\l+(l
-1)\a_1+l\a_2}.$$ Let
$$y_1=\frac12[[f_2^{(l-1)},f_1^{(l)}],f_1^{(l)}]f_1^{(l-3)}.$$
Then $y_1f_2^{(2l-2)}f_1^{(l-1)}=x_{l-1,l-1}$. Using 1.1 (c) we
see that $y_1f_2^{(l-2)}\tol$ is a primitive element of weight
$\l-(3l -3)\a_1-(2l-3)\a_2$. Note that $y_1f_2^{(l-2)}=
{x_{l-1,l-1}}/{f_2^{(l)}f_1^{(l-1)}}\in\uxin.$

\def\tbu{\tilde{\bold u}}

\medskip

 \noindent(2)
Now we consider the homomorphism: $$\tilde Z(\l)\to\tilde
Z(\l+2\a_2),\quad \tol\to t_2=f_2^{(2)}\tilde 1_{\l+2\a_2}.$$

Let
$$x_2=(f_2^{(l-1)}f_1^{(l)}-f_1^{(l)}f_2^{(l-1)})f_1\in\uxin.$$
  Using 1.1 (c) we see that
$x_2\tol$ is a  primitive element of weight $\l-(l +1)\a_1-
(l-1)\a_2$.

 Let $y_2=f_1^{(l-3)}x_2\in\uxin.$
 According to 1.1 (c) we know that $y_2\tol$ is a primitive
 element of weight $\l-(2l-2)\a_1-(l -1)\a_2$. It is easy to see
 that $y_2\tol$ is maximal.

\medskip

\noindent(3) Now we consider the homomorphism: $$\tilde
Z(\l)\to\tilde Z(\l+3\a_1+2\a_2),$$ $$ \tol\to
t_3=f_2^{(2)}f_1^{(3)}\tilde 1_{\l+3\a_1+2\a_2}.$$ Let
$x_3=\frac{x_{3,l-1}}{f_2^{(2)}f_1^{(3)}}\in\uxin.$ Using 1.1 (c)
we know that $x_3\tol$ is a primitive element of weight $\l-(2l
+1)\a_1-(2l -1)\a_2$.

\medskip

Consider the homomorphism: $$\tilde Z(\l)\to\tbu f_2^{(l+2)}f_1^{(3)}\tilde 1_{\l+3\a_1+(l+2)\a_2},$$ $$ \tol\to
t_3=f_2^{(l+2)}f_1^{(3)}\tilde 1_{\l+3\a_1+(l+2)\a_2}.$$ Let
$y_3=[[f_2^{(l-1)},f_1^{(l)}],f_1^{(l)}]f_1=\frac{x_{3,l-1}}{f_2^{(l+2)}f_1^{(3)}}\in\uxin.$
Using 1.1 (c) we know that $y_3\tol$ is a primitive element of
weight $\l-(2l +1)\a_1-(l -1)\a_2$.

\medskip

\noindent(4) We consider the homomorphism: $$\tbu f_1\tol\to\tbu
f_2f_1^{(3)}\tilde 1_{\l+2\a_1+\a_2} ,$$ $$ f_1\tol\to t_4=f_2
f_1^{(3)}\tilde 1_{\l+2\a_1+\a_2}.$$ Let $x_4=
f_1^{(l-1)}f_2^{(l)}f_1\in\uxin.$ Using 1.1 (c) we know that
$x_4\tol$ is a primitive element of weight $\l-l\a_1-l\a_2$.

Let $y_4=\frac{x_{3,l-2}}{f_2f_1^{(3)}}\in\uxin.$  Using 1.1
(c) we know that $y_4f_1\tol$ is a primitive element of weight
$\l-2l\a_1-(2l-2)\a_2$.

\medskip

 \noindent(5)
Now we consider the homomorphism: $$\tilde Z(\l)\to\tilde
Z(\l+(l+2)\a_1+(l+1)\a_2),\quad \tol\to
t_3=f_1^{(l-1)}f_2^{(l+1)}f_1^{(3)} \tilde
1_{\l+(l+2)\a_1+(l+1)\a_2}.$$ Let
$$x_5=(f_2^{(l-2)}f_1^{(l)}-f_1^{(l)}f_2^{(l-2)}) \in\uxin.$$
  Using 1.1 (c) we see that
$x_5\tol$ is a  primitive element of weight $\l-l\a_1- (l-2)\a_2$.

\medskip

\noindent(6) Consider the homomorphism $$\tbu f_1 \tol\to \tbu
f_2f_1^{(3)}f_2^{(2)} \tilde 1_{\l+2\a_1+3\a_2}$$ $$ f_1 \tol\to
f_2f_1^{(3)}f_2^{(2)} \tilde 1_{\l+2\a_1+3\a_2}.$$ Let
$x_6=\frac{x_{l-1,2}}{f_2f_1^{(3)}f_2^{(2)}}\in\uxin$ By 1.1 (c),
$x_6f_1\tol$ is primitive and is of weight $\l-2l\a_1-l\a_2$.

\bigskip

Note that  the element $m=f_2^{(l
 -2)} f_1^{(2l-3)}f_2^{(l -1)}f_1 \tol$  generates the
unique irreducible submodule of $\tzl$. By comparing
the weights of the   following 6 elements with the weight of
$m$,
$$\frac{x_{l-1,l-1}}{f_1^{(l-1)}}\tol,\quad\frac{x_{l-1,l-1}}{f_2^{(l)}
f_1^{(l-1)}}\tol,\qquad
\frac{x_{3,l-1}}{f_2^{(2)}f_1^{(3)}}\tol, $$ $$\frac{x_{3,l-1}}{f_2^{(l+2)}
f_1^{(3)}}\tol,\qquad
f_2^{(l-2)}f_1^{(l-1)}f_2^{(l)}f_1\tol,
\quad\frac{x_{l-1,2}}{f_2
f_1^{(3)}f_2^{(2)}}f_1\tol,$$ we see
that there are no maximal elements in $\tzl$ that have the same
weight with any of  above 6 elements.

Now we show that there are no maximal elements in $\tzl$ that have the same
weight with any of  other 3 elements in (ii) by assuming (iii).

By (iii), $\tzl$ has only one composition factor isomorphic to
$\tilde L(\l-(l+1)\a_1-(l-1)\a_2)$. Suppose that there is a maximal element
$m$
in $\tzl$ of weight $\l-(l+1)\a_1-(l-1)\a_2$. Then $m$ is in $\tbu
y$, here $y=[f_2^{(l-1)},f_1^{(l)}]f_1\tol$. It is clear that
$\tbu y\subset\tbu^-y+\tbu^-f_1^{(l  -3)}f_2^{(l -1)}f_1 \tol$.
Thus $m=ay+bf_1^{(3)}f_1^{(l  -3)}f_2^{(l -1)}f_1 \tol$ for some
$a,b$ in ${\bold Q}(\xi)$. Thus $m=ay$. But $y$ is not maximal. So
there are no maximal elements in $\tzl$ that have the same
weight with $[f_2^{(l-1)},f_1^{(l)}]f_1\tol$.

Similarly, we see that there are no maximal elements in $\tzl$ that have the same
weight with any of  $[f_2^{(l-2)},f_1^{(l)}]\tol,\
f_1^{(l-1)}f_2^{(l)}f_1\tol$.

\noindent(iii) Using (i), (ii) and  1.1 (b), we see that $\tzl$
has the 20 composition factors. By a comparison of dimensions we
know that $\tzl$ has only the 20 composition factors.

 The theorem is proved.
\bigskip

\noindent {\bf Theorem 3.3.} {\sl   Let $a,b$ be integers and
$\l=(la+l-2 ,lb+1) $. Then}

\noindent(i){\sl  The following  elements are maximal     in
$\tilde Z(\l)$: $$\tol,\qquad f_1^{(l-1)} \tol,\qquad
f_2^{(2)}\tol,\qquad f_2f_1^{(l -1)} \tol,$$ $$
f_1^{(3)}f_2^{(2)}\tol,\qquad
 f_1^{(3)}f_2^{(l+1)}f_1^{(l-1)} \tol.$$
  $$
f_2 f_1^{( 3)}f_2^{( 2)}\tol,\qquad f_2^{(2)} f_1^{(l+3)}f_2^{(l
+1)}f_1^{(l-1)} \tol,$$} $$\frac{f_1^{( 3)}f_2 }{f_1}\tol,\qquad
\frac{f_2^{(2)}f_1^{(3)}f_2 }{f_1}\tol,\qquad
\frac{x_{3,l-2}}{f_2^{(l-2)} }\tol .$$

\noindent(ii){\sl  The following   elements are  primitive
elements in $\tilde Z(\l)$ but not maximal:
$$[f_2^{(2)},f_1^{(l)}]\tol,\quad
\frac{x_{l-1,l-1}}{f_2^{(l-2)}f_1^{(2l-3)}f_2^{(l-1)}}\tol, \quad
\frac{f_2^{(l-1)}f_1^{(l)}f_2}{f_1}\tol,\quad
\frac{f_1f_2^{(l-1)}f_1^{(l)}f_2}{f_1}\tol,$$
$$\frac{[f_2^{(2)},f_1^{(l)}]f_1^{(3)}f_2  }{f_1}\tol,
\quad\frac{x_{3,l-1}}{f_1f_2^{(l-1)}  }\tol,\quad
\quad\frac{x_{3,l-1}}{f_1^{(l+1)}f_2^{(l-1)} }\tol,\quad
\frac{x_{3,l-1}}{f_1^{(2l+1)}f_2^{(l-1)}  }\tol,$$
$$\frac{x_{l-1,l-1}}{
f_1^{(l-3)}f_2^{(l-1)}}f_2^{(2)}\tol .$$

Moreover there are no maximal     elements in $\tilde Z(\l)$ which
have the same weight with any of above 9 elements.}

\def\tl{\tilde L}

\noindent(iii) {\sl The maximal and primitive          elements in
(i-ii) provide 20 composition factors of $\tilde Z(\l)$, which are
$$\tilde L(\l),\qquad\tilde L(\l- (l-1)\a_1),\qquad\tilde L(\l-2
\a_2),$$ $$\tilde L(\l- (l-1)\a_1- \a_2),\qquad\tilde L(\l-3\a_1-2
\a_2)$$
 $$\tilde L(\l- (l+2)\a_1-(l+1)\a_2),\qquad\tilde L(\l-3\a_1-
3\a_2),$$ $$\tilde L(\l-(2l+2)\a_1-(l+3)\a_2)\qquad\tilde
L(\l-2\a_1- \a_2),$$ $$
\tl(\l-2\a_1-3\a_2),\qquad\tl(\l-l\a_1-2\a_2),
\qquad\tl(\l-(2l-1)\a_1-l\a_2), $$
$$\tl(\l-(l-1)\a_1-l\a_2),\qquad\tl(\l-l\a_1-l\a_2),
\qquad\tl(\l-(l+2)\a_1 -3\a_2),$$ $$
 \tl(\l-(2l+3)\a_1-(l+2)\a_2),\qquad\tl(\l-(l+3)\a_1-(l+2)\a_2),
 \qquad\tl(\l-3
 \a_1-(l+2)\a_2),$$ $$
 \tl(\l-(2l+2)\a_1-(l+1)\a_2),\qquad\tl(\l-(3l-1)\a_1-2l\a_2) .$$
 Moreover, $\tilde Z(\l)$ has
only the 20 composition factors.}

\def\g{\gamma}

\noindent{\it Proof:} (i) According to 1.1 (e-f), we see that the first
8 elements in (i) are maximal.

Consider the homomorphism: $$\tilde Z(\l)\to\tilde
Z(\l+\a_1),\quad \tol\to t_1=f_1\tilde 1_{\l+\a_1}.$$

Since $f_1^{( 3)}f_2 $ is in $\bu^-f_1$, using 1.1 (c) we see that
$\frac{f_1^{( 3)}f_2 }{f_1}\tol$ and   $\frac{f_2^{(2)}f_1^{(
3)}f_2 }{f_1}\tol$ are primitive elements of weights $\l-2\a_1-
\a_2$ and $\l-2\a_1-3\a_2$ respectively.

\medskip

Now we consider the homomorphism: $$\tilde Z(\l)\to\tilde
Z(\l+(l-2)\a_2)),\quad \tol\to t_2=f_2^{(l-2)}\tilde
1_{\l+(l-2)\a_2}.$$ Using 1.1 (c) we see that
$f_1^{(3)}f_2[f_2^{(l)},f_1^{(l)}]f_1^{(l-1)}\tol=\frac{x_{3,l-2}}{f_2^{(l-2)}}\tol$ is a primitive element of
weight $\l-(2l+2)\a_1-(l+1)\a_2$.

\medskip

One can check directly
that the three elements are maximal.

\medskip

\noindent(ii) Now we argue for (ii).

 \noindent(1)
Consider the homomorphism: $$\tilde Z(\l)\to\tilde Z(\l+
\a_1),\quad \tol\to t_1=f_1 \tilde 1_{\l+ \a_1}.$$ Using Theorem
3.1 and 1.1 (c) we see that the first 8 elements are primitive.

\medskip

\noindent(2) Now we consider the homomorphism: $$\tbu
f_2^{(2)}\tol\to\tbu f_1^{(l-3)}f_2^{(l-1)}\tilde
1_{\l+(l-3)\a_1+(l-3)\a_2} ,$$ $$ f_2^{(2)}\tol\to
f_1^{(l-3)}f_2^{(l-1)}\tilde 1_{\l+(l-3)\a_1+(l-3)\a_2} .$$  Using
1.1 (c) we know that
$\frac{x_{l-1,l-1}}{f_1^{(l-3)}f_2^{(l-1)}}f_2^{(2)}\tol$ is a
primitive element of weight $\l-(3l-1)\a_1-2l\a_2$.

\bigskip
Consider the homomorphism: $$\tilde Z(\l)\to\tilde Z(\l+
\a_1)),\quad \tol\to t_1=f_1 \tilde 1_{\l+ \a_1}.$$ It is easy to
see that all the 9 primitive elements have non-zero image. By Theorem 3.1  (ii) and 1.1 (e) we know that there are
no maximal
elements in $\tilde Z(\l)$ which have the same weight with any of
the 9 primitive elements.

\noindent(iii) Using (i), (ii) and  1.1 (b), we see that $\tzl$
has the 20 composition factors. By a comparison of dimensions we
know that $\tzl$ has only the 20 composition factors.

The theorem is proved.
\bigskip

\noindent {\bf Theorem 3.4.} {\sl   Let $a,b$ be integers and
$\l=(la+2 ,lb+l-2) $. Then}

\noindent(i){\sl  The following  elements are maximal     in
$\tilde Z(\l)$: $$\tol,\qquad f_1^{(3)} \tol,\qquad
f_2^{(l-1)}\tol,\qquad f_2^{(2)}f_1^{(3)} \tol,$$ $$ f_1
f_2^{(l-1)}\tol,\qquad
 f_1 f_2^{(2)}f_1^{(3)} \tol.$$
  $$
f_2^{(2)} f_1^{( l+1)}f_2^{(l-1)}\tol,\qquad f_2^{(l-1)}
f_1^{(2l+1)}f_2^{(l +2)}f_1^{(3)} \tol,$$} $$\frac{f_2^{(2)}f_1
}{f_2}\tol,\qquad \frac{f_1^{(3)}f_2^{(2)}f_1 }{f_2}\tol.$$

\noindent(ii){\sl  The following   elements are  primitive
elements in $\tilde Z(\l)$ but not maximal:
$$[f_2^{(l-1)},f_1^{(l)}]\tol,\quad
\frac{x_{l-1,l-1}}{f_1^{(l-3)}f_2^{(l-2)}f_1^{(l-1)}}\tol, \quad
\frac{x_{3,l-2}}{f_2f_1^{(l-1)}f_2^{(l-2)}}\tol,\quad
[f_1^{(3)},f_2^{(l)}]\tol,$$ $$\frac{x_{3,l-1}
}{f_2f_1^{(2l)}f_2^{(l-1)}}\tol,
\quad\frac{x_{3,l-1}}{f_2f_1^{(l)}f_2^{(l-1)} }\tol,\quad \quad
f_1[f_2^{(l-1)},f_1^{(l)}] \tol,\quad
\frac{x_{l-1,2}}{f_2f_1^{(l-1)} }\tol,$$
$$[f_2^{(2)},f_1^{(l)}]f_1^{(3)}\tol,\quad\frac{x_{l-1,l-1}}{
f_1^{(l-3)}f_2^{(2l-2)}f_1^{(l-1)}}\tol .$$

Moreover no maximal     element in $\tilde Z(\l)$ has the same
weight as any of above 10 elements.}

\def\tl{\tilde L}

\noindent(iii) {\sl The maximal and primitive          elements in
(i-ii) provide 20 composition factors of $\tilde Z(\l)$, which are
$$\tilde L(\l),\qquad\tilde L(\l- 3\a_1),\qquad\tilde L(\l-(l-1)
\a_2),$$ $$\tilde L(\l- 3\a_1-2 \a_2),\qquad\tilde L(\l-
\a_1-(l-1) \a_2)$$
 $$\tilde L(\l- 4\a_1-2\a_2),\qquad\tilde L(\l-(l+1)\a_1-
(l+1)\a_2),\qquad\tilde L(\l-(2l+4)\a_1-(2l+1) \a_2),$$ $$
\tl(\l-\a_1- \a_2),\qquad\tl(\l-4\a_1-\a_2),
\qquad\tl(\l-l\a_1-(l-1)\a_2), $$
$$\tl(\l-2l\a_1-(2l-1)\a_2),\qquad\tl(\l-(l+3)\a_1-l\a_2),
\qquad\tl(\l-3\a_1 -l\a_2),$$ $$
 \tl(\l-4\a_1-(l+1)\a_2),\qquad\tl(\l-(l+4)\a_1-(l+1)\a_2),
 \qquad\tl(\l-(l+1)
 \a_1-(l-1)\a_2),$$ $$
 \tl(\l-(l+3)\a_1-(l+2)\a_2),\qquad\tl(\l-(l+3)\a_1-2\a_2),
 \qquad \tl(\l-2l\a_1-(l-1)\a_2).$$
 Moreover, $\tilde Z(\l)$ has
only the 20 composition factors.}

\def\g{\gamma}

\noindent{\it Proof:} (i) According to 1.1 (e-f), we see the first
8 elements in (i) are maximal.

Consider the homomorphism: $$\tilde Z(\l)\to\tilde
Z(\l+\a_2),\quad \tol\to t_1=f_2\tilde 1_{\l+\a_2}.$$

Since $f_2^{(2)} f_1 $ is in $\bu^-f_2$, using 1.1 (c) we see that
$\frac{f_2^{(2)} f_1  }{f_2}\tol$ and   $\frac{ f_1^{(
3)}f_2^{(2)} f_1  }{f_2}\tol$ are primitive elements of weights
$\l-\a_1- \a_2$ and $\l-4\a_1-\a_2$ respectively. One can check
directly that the two elements are maximal.

\noindent(ii) Now we argue for (ii).

 \noindent(1)
Consider the homomorphism: $$\tilde Z(\l)\to\tilde Z(\l+
\a_2)),\quad \tol\to t_1=f_2 \tilde 1_{\l+ \a_2}.$$ Using Theorem
3.1 and 1.1 (c) we see that the first 9 elements are primitive.

\medskip

 \noindent(2)
 Now we consider the homomorphism: $$\tzl\to\tbu f_1^{(l-3)}f_2^{(2l-2)}f_1^{(l-1)}
 \tilde
1_{\l+(2l-4)\a_1+(2l-2)\a_2} ,$$ $$  \tol\to
f_1^{(l-3)}f_2^{(2l-2)}f_1^{(l-1)}
 \tilde
1_{\l+(2l-4)\a_1+(2l-2)\a_2} .$$  Using 1.1 (c) we know that
$[[f_2^{(l-1)},f_1^{(l)}],f_1^{(l)}]\tol=\frac{x_{l-1,l-1}}{f_1^{(l-3)}f_2^{(2l-2)}f_1^{(l-1)}}\tol$
is a primitive element of weight $\l-2l\a_1-(l-1)\a_2$.

Consider the homomorphism: $$\tilde Z(\l)\to\tilde Z(\l+
\a_2)),\quad \tol\to t_1=f_2 \tilde 1_{\l+ \a_2}.$$ It is easy to
see that all the 10 primitive elements have non-zero image. By Theorem 3.1  (ii) and 1.1 (e) we know that there are
no maximal
elements in $\tilde Z(\l)$ which have the same weight with any of
the 10 primitive elements.

\bigskip

\noindent(iii) Using (i), (ii) and  1.1 (b), we see that $\tzl$
has the 20 composition factors. By a comparison of dimensions we
know that $\tzl$ has only the 20 composition factors.

The theorem is proved.

\bigskip
\section{ Weyl Modules for Type $B_2$}
\bigskip
\def\bzpn{\bold Z^2_+}
 For $\l=(\l_1,\l_2)\in\bzpn$ we denote by $I\ll$
the left ideal of $\Uxi$ generated by all $\ea\ (a\geq 1) $,
$\f^{(a_i)}\ (a_i\geq \l_i+1)$,
$\k-\xi^{\l_i},\kca-{\left[\begin{array}{r} \l_i+c\\ a\
\
\
\end{array}\right] }_{\xi^i}$. The Weyl
module $V(\l)$ of $\Uxi$ is defined to be $\Uxi/I\ll$, its
dimension is $(\l_1+1)(\l_2+1)(\l_1+\l_2+2)(\l_1+2\l_2+3)/6$. Let
$v\ll$ be a nonzero element in $V(\l)\ll$.
 We can work out the maximal and primitive elements in $V(\l)$ as
 in section 3 (cf. [X2]). We omit the
 results here.

\bigskip
\noindent{\bf Acknowledgement:} \rm The work was completed during
my visit to the University of Sydney. It is a great pleasure to
thank Professor G. Lehrer for the invitation.  Part of the work
was done   during my visit to Bielefeld University. I am grateful
to the SFB 343 in Bielefeld University for financial support.

\bigskip\bigskip\bigskip\bigskip

\bibliographystyle{amsalpha}
{}
\bigskip
\bigskip\bigskip

\end{document}